\documentclass[12pt]{article}
\usepackage{mathptmx}
\usepackage{amssymb,latexsym,amsmath,amsthm}
\usepackage{stmaryrd}
\usepackage{setspace}
\usepackage{bm}
\usepackage{bbm}
\usepackage{dsfont}
\usepackage{graphicx}
\usepackage{epsfig}
\usepackage{tikz}
\usetikzlibrary{shapes,backgrounds,calc}
\usepackage{cite}
\usepackage{enumerate}
\usepackage{hyperref}
\usepackage[top=1in, bottom=1in, left=1in, right=1in]{geometry}
\usepackage{boxedminipage}
\usepackage{verbatim}
\usepackage{enumitem}
\usepackage{fancyhdr}
\usepackage[utf8]{inputenc}
\usepackage{color}
\usepackage{mathtools}
\definecolor{myGreen}{rgb}{0,.5,0}
\definecolor{black}{rgb}{0,0,0}
\definecolor{myYellow}{rgb}{.6,.6,.1}
\definecolor{Cyan}{rgb}{.2,.6,.7}
\definecolor{white}{rgb}{1,1,1}
\definecolor{Purple}{rgb}{.4,0,1}
\definecolor{deepblack}{rgb}{.53,.29,.24}
\definecolor{Black}{rgb}{0,0,0}
\definecolor{Grey}{rgb}{.45,.45,.45}

\newcommand{\noi}{\noindent}
\newcommand{\op}[1]{\left(#1\right)}

\newcommand{\abs}[1]{\left\vert#1\right\vert}

\newcommand{\N}{\mathbb{N}}
\newcommand{\R}{\mathbb{R}}

\newcommand{\ds}{\displaystyle}

\newcommand{\multiset}[2]{\ensuremath{\left(\kern-.3em\left(\genfrac{}{}{0pt}{}{#1}{#2}\right)\kern-.3em\right)}}
\newcommand{\xmin}{x_\text{min}}
\newcommand{\xmax}{x_\text{max}}
\definecolor{cb1}{RGB}{230,159,0}
\definecolor{cb4}{RGB}{86,180,233}
\definecolor{cb2}{RGB}{0,158,115}
\definecolor{cb5}{RGB}{0,114,178}
\definecolor{cb3}{RGB}{213,94,0}
\definecolor{cb6}{RGB}{204,121,167}
\tikzset{vtx/.style={circle,black,fill=black,inner sep=1pt,minimum size=7pt}}
\tikzset{whiteVtx/.style={circle,white,fill=white,inner sep=1pt}}
\tikzset{goldVtx/.style={circle,white,fill=CUGold,inner sep=1pt}}
\tikzset{cb1Vtx/.style={circle,white,fill=cb1,inner sep=1pt}}
\tikzset{cb2Vtx/.style={circle,white,fill=cb2, inner sep=1pt}}
\tikzset{cb3Vtx/.style={circle,white,fill=cb3, inner sep=1pt}}
\tikzset{cb4Vtx/.style={circle,white,fill=cb4, inner sep=1pt}}
\tikzset{cb5Vtx/.style={circle,white,fill=cb5, inner sep=1pt}}
\tikzset{blackVtx/.style={circle,white,fill=black,inner sep=1pt}}
\tikzset{type/.style={draw=cb6,very thick,fill=yellow}}
\tikzset{type2/.style={draw=cb6,very thick,fill=yellow,inner sep=1pt,minimum size=8pt}}
\tikzset{type3/.style={circle,draw=black,fill=white,inner sep=1pt,minimum size=7pt}}
\tikzset{e1/.style={line width=4pt,cb3}}
\tikzset{e2/.style={line width=4pt,cb4}}
\tikzset{e3/.style={line width=4pt,cb2}}
\tikzset{eq/.style={line width=4pt,gray!40}}
\tikzset{dedge/.style={->,line width=1.3pt}} % directed edge
\newtheorem{theorem}{Theorem}
\newtheorem{thm}{}[section]
\newtheorem{lemma}[thm]{Lemma}
\newtheorem{claim}[thm]{Claim}
\newtheorem{conjecture}[theorem]{Conjecture}

\newtheorem{prop}[thm]{Proposition}

\newtheorem{case}{Case}

\numberwithin{subcase}{case}
\theoremstyle{definition}

\graphicspath{{./graphics/}}

\begin{document}

\title{Inducibility of the Net Graph}

\author{Adam Blumenthal\thanks{Department of Mathematics and Computer Science, Westminster College, E-mail: {\tt blumenam@westminster.edu} .}
\and
Michael Phillips\thanks{Department of Mathematical and Statistical Sciences, University of Colorado Denver, E-mail: {\tt Michael.2.Phillips@ucdenver.edu} .} }

\maketitle

\begin{abstract}
A graph $F$ is called a fractalizer if for all $n$ the only graphs which maximize the number of induced copies of $F$ on $n$ vertices are the balanced iterated blow ups of $F$. 
While the net graph is not a fractalizer, we show that the net is nearly a fractalizer.
Let $N(n)$ be the maximum number of induced copies of the net graph among all graphs on $n$ vertices. 
For sufficiently large $n$ we show that, $N(n) = x_1\cdot x_2 \cdot x_3 \cdot x_4 \cdot x_5 \cdot x_6 + N(x_1) + N(x_2) + N(x_3) + N(x_4) + N(x_5) + N(x_6)$ where $\sigma x_i = n$ and all $x_i$ are as equal as possible.
Furthermore, we show that the unique graph which maximizes $N(6^k)$ is the balanced iterated blow up of the net for $k$ sufficiently large.
We expand on the standard flag algebra and stability techniques through more careful counting and numerical optimization techniques.
\end{abstract}

\section{Introduction}

In 1975, Pippenger and Golumbic \cite{pipGol} proposed the \emph{inducibility problem}: determine the maximum possible density of induced copies of a $k$-vertex graph $H$ that can be contained in an $n$-vertex graph. Here, we take the density to be the number of copies of $H$ divided by $\binom{n}{k}$, the number of $k$-vertex induced subgraphs of the host graph. For a $k$-vertex graph $H$ and $n$-vertex graph $G$, we denote by $I(H,G)$ the number of induced copies of $H$ in $G$, and by $i(H,G) := I(H,G)/\binom{n}{k}$ the density of $H$ in $G$. If $G$ is an $n$-vertex graph which maximizes $I(H,G)$ over all $n$-vertex graphs, then we say that $G$ is extremal and introduce the notations $I(H,n) = I(H,G)$ and $i(H,n) = i(H,G)$. Finally, we define the inducibility of $H$ as $\ds\lim_{n\to\infty} \frac{I(H,n)}{\binom{n}{k}} = \lim_{n\to\infty} i(H,n)$. 

Pippenger and Golumbic \cite{pipGol} proved for all $k$-vertex graphs $H$ that $i(H) \geq \frac{k!}{k^k-k} = (1 + o(1))\frac{k!}{k^k}$ by construction. Although they use different verbiage, the authors use an \emph{iterated, balanced blow-up} of a graph $H$ as their construction for the lower-bound on $i(H)$. A \emph{blow-up} of a graph $H = (V(H)=\{v_1,\dots,v_k\},E(H))$ by replacing each vertex $v_i$ with a graph $H_i$ so that for any distinct $i,j \in [k]$, a pair of vertices $x \in H_i$ and $y \in H_j$ are adjacent if and only if $v_iv_j \in E(H)$. We say that a blow-up is \emph{balanced} (or nearly balanced) if $\abs{|H_i| - |H_j|} \leq 1$ for all $i,j \in [k]$, and call a blow-up \emph{iterated} if each $H_i$ is itself a blow-up of $H$, each part in each $H_i$ is a blow-up, etc.

Until Razborov \cite{razborov} introduced the flag algebra method in 2007, little progress had been made toward the resolution of the inducibility problem. Since then, the problem has been closed for several small graphs \cite{balHuLidPfen,hirst,EvenZohar2014}. Fox, Huang and Lee \cite{foxhuanglee}, and separately Yuster \cite{yuster} resolved this problem for almost all graphs by considering $G_{n,p}$ for arbitrary fixed $p$ as $n \to \infty$. Here, $G_{n,p}$ is a vector of $\binom{n}{2}$ Bernoulli($p$) random variables, each determining the presence or absence of an edge in an $n$-vertex graph. In order to fully express their result, we also need the following definition.

A graph $F$ is a \textbf{fractalizer} if, for each $n$, all graphs on $n$ vertices which maximize the number of induced copies of $F$ are balanced iterated blow-ups of $F$. The result by Fox, Huang and Lee states that $G_{n,p}$ is a fractalizer (or fractalizes) with high probability. However, other than $K_n$ and $\overline{K_n}$, no small or explicit, large graphs are known to be fractalizers. One issue is with the strength of the definition. It has been shown by Lidick\'{y}, Mattes, and Pfender \cite{c5frac} that $C_5$ is \emph{almost} a fractalizer: \begin{itemize}
	\item for $n$ sufficiently large, all extremal constructions are balanced $C_5$ blow-ups,
	\item balanced iterated blow-ups are extremal for all $n$, but
	\item there is a structure on 8 vertices (the M\"{o}bius ladder) which has the same number of induced copies of $C_5$ as its balanced blow-up.
\end{itemize}

\begin{figure}[t] \centering
        \newcommand{\mininet}{
        \begin{tikzpicture}[rotate = -90]
        \foreach \j in {1,2,3}{
            \node[scale = .25, fill = black, circle] at ({360/3 * \j}:.15) (II\j) {};
            \node[scale = .25, fill = black, circle] at ({360/3 * \j}:.35) (OO\j) {};
            \draw[thick] (II\j) -- (OO\j);
        }
        \draw[thick](II1)--(II2)--(II3)--(II1);
        \end{tikzpicture}
        }
        
        \begin{tikzpicture}
        
        \begin{scope}[xshift = -1in, rotate = -90]
        \foreach \i in {1,2,3}{
            \node[scale = .5, fill = black, circle] at ({360/3 * \i}:.75) (I\i) {};
            \node[scale = .5, fill = black, circle] at ({360/3 * \i}:1.5) (O\i) {};
            \draw[thick] (I\i) -- (O\i);
        }
        \draw[thick] (I1)--(I2)--(I3)--(I1);
        \end{scope}
        
        \begin{scope}[xshift = 1in, rotate = -90]
        \foreach \i in {1,2,3}{
            \node[scale = .75, draw, circle] at ({360/3 * \i}:.75) (I\i) {\mininet};
            \node[scale = .75, draw, circle] at ({360/3 * \i}:2) (O\i) {\mininet};
            \draw[line width = .3cm] (I\i) -- (O\i);
            }
        \draw[line width = .3cm] (I1)--(I2)--(I3)--(I1);
        \end{scope}
        \end{tikzpicture}
\caption{The net graph and the iterated blow-up of the net graph on 36 vertices.}\label{netGraph}
\end{figure}
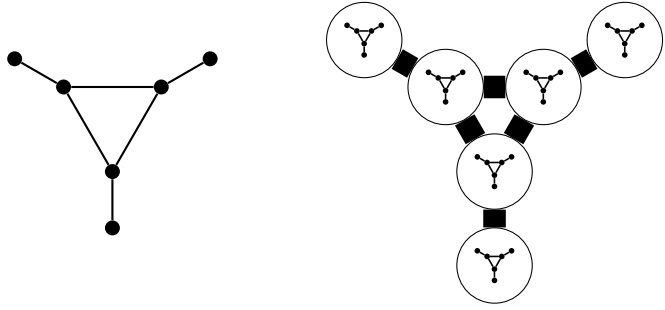

In this paper, we determine the inducibility of the net graph (as seen in Figure \ref{netGraph}), as well as the unique $n$-vertex graphs which maximize the density of induced net graphs, or \emph{nets}, when $n$ is a power of 6. Our work is heavily influenced by \cite{balHuLidPfen} wherein the Balough, Hu, Lidick\'{y}, and Pfender resolve one case of another 1975 question of Pippenger and Golumbic \cite{pipGol} regarding the inducibility of cycle graphs. In this paper, we build on their methods to include careful analysis of particular subgraph densities to determine for sufficiently large $n$ the graph which maximizes the density of the net graph.

We show that the net graph satisfies a weaker condition (see \cite{c5frac} for more discussion regarding this topic): the balanced iterated blow-up of the net uniquely maximizes the density of the net graph when $n$ is a power of 6, and for $n$ large enough, the graphs which maximize the density of the net are balanced, iterated blow-ups of the net. We note that the graph obtained by adding a pendant to each vertex in a $K_4$ has the same number of induced copies of the net graph as a balanced iterated blow-up of the net on 8 vertices, so the net cannot be a fractalizer. As it stands, there are still no known nontrivial fractalizers.

Inducibility has natural extensions to directed graphs as well. Falgas-Ravry and Vaughan \cite{MR2988862} considered inducibility of small outstars, with an extension to all outstars by Huang in \cite{Huang2014} before being further generalized to other types of stars by Hu, Ma, Norin, and Wu \cite{HMNW}. Orientations of short paths and the $C_4$ are explored in \cite{choi} and \cite{huVolec}, respectively. The inducibility of directed graphs on at most 4 vertices was explored in \cite{bozyk2020inducibility}, and the corresponding tournaments were completely resolved in \cite{myTournaments}. Other recent work includes closing bipartite graphs on 5 vertices in \cite{k2111} and an exploration of inducibility for $d$-ary trees in \cite{dary}.

Additionally, we would like to provide the reader with a better sense of how one would approach problems like those presented in this paper. We believe that many steps of this process can (and should) be stream-lined, either by explicitly providing more generalized results or proposing stages which can be automated. To that end, we state and prove an important property of vertices in extremal constructions of inducibility problems: every vertex in an extremal construction should be in roughly the same number of induced copies of the target structure.

\begin{lemma}\label{AvgArg}
Let $H$ be a $k$-vertex graph, and let $G$ be an $n$-vertex graph which satisfies $I(H,G) = I(H,n)$. Every $v \in V(G)$ is contained in at least $(L + o(1))\binom{n}{k-1}$ copies of $H$ where $L$ is any lower-bound on $i(H)$.
\end{lemma}

\noindent
\begin{proof}
Let $L$ be any lower-bound for the inducibility of the $k$-vertex graph $H$. If we denote by $H^u$ the number of nets containing a vertex $u$ in an extremal graph $G$, then $\sum_{u \in V(G)} H^u \geq k \cdot (L + o(1)) \binom{n}{k}$, implying that the average $\overline{H^u}$ of $H_u$ over $V(G)$ is at least $(L + o(1))\binom{n}{k-1}$.

Let $u,v \in V(G)$ be given arbitrarily, and denote by $H^{uv}$ the number of nets containing both of $u$ and $v$. Trivially, we have that $H^{uv} \leq \binom{n-2}{k-2}$. Construct $G'$ from $G$ by deleting $v$ and duplicating $u$ as $u'$. As $G$ is extremal, we have the following sequence of inequalities:
$$
0 \geq I(H,G') - I(H,G)
	\geq H^u - H^v - H^{uv}
	\geq H^u - H^v - \binom{n-2}{k-2}.
$$Therefore, we have for every pair $u,v \in V(G)$ that $|H^u - H^v| \leq \binom{n-2}{k-2}$. As some vertex must be in at least $\overline{H^u}$ induced copies of $H$, it follows that $H^u$ is at least
\begin{align*}
\overline{H^u} - \binom{n-2}{k-2}
	&\geq (L + o(1))\binom{n}{k-1} - \binom{n-2}{k-2} \\
	&= (L + o(1))\binom{n}{k-1} - o(n^{k-1}) \\
	&= (L + o(1))\binom{n}{k-1}
\end{align*}for every $u \in V(G)$, as desired.
\end{proof}

The main result of this paper is Theorem \ref{net1}, which follows from Theorem \ref{net2} (below) using a standard argument that is largely reproducible for other graphs.

\begin{theorem}\label{net1}
For $k \geq 1$, the unique graph on $6^k$ vertices which maximizes the number of induced copies of the net graph is a balanced, iterated blow-up of the net.
\end{theorem}

In order to prove Theorem \ref{net1}, we first prove Theorem \ref{net2}, which is sufficient for determining the unique limit object maximizing the density of induced copies of the net graph. A proof sketch for Theorem \ref{net2} follows in Section \ref{ProofSketch}, and the proof of Theorem \ref{net2} fills most of Section \ref{ProofNotSketch} with the final arguments to conclude Theorem \ref{net1} at the very end.

\begin{theorem}\label{net2}
There exists $n_0$ such that for every $n \geq n_0$
$$
N(n) = x_1\cdot x_2\cdot x_3\cdot x_4\cdot x_5\cdot x_6 + \sum_{1 \leq i \leq 6} N(x_i)
$$where $\sum x_i = n$ and $|x_i - x_j| \in \{0,1\}$ for all $i,j \in [6]$. Moreover, if $G$ maximizes the density of nets among $n$-vertex graphs, then $G$ is a balanced blow-up of the net.
\end{theorem}

\section{Proof Sketch for Theorem 2}\label{ProofSketch}

We begin the proof of Theorem \ref{net2} by showing that the extremal construction is in some sense very close to a balanced blow-up of the net; we would need only to delete less than $0.2\%$ of the vertices and add/remove at most $0.001\%$ of the edges in an extremal graph for it to become a blow-up of the net with part sizes \emph{almost} balanced (Claim \ref{claim4}). Once this has been accomplished, we argue that actually no edges should be added or removed (Claim \ref{claim5}), and then further that no vertices need to be deleted (Claims \ref{claim6} through \ref{claim8}), implying that the extremal construction is precisely a blow-up of the net. Once we have this result, we argue that the parts are all asymptotically the same size using basic analytic techniques, then that they differ by at most one vertex (Claim \ref{claim9}), concluding the proof of Theorem 2.

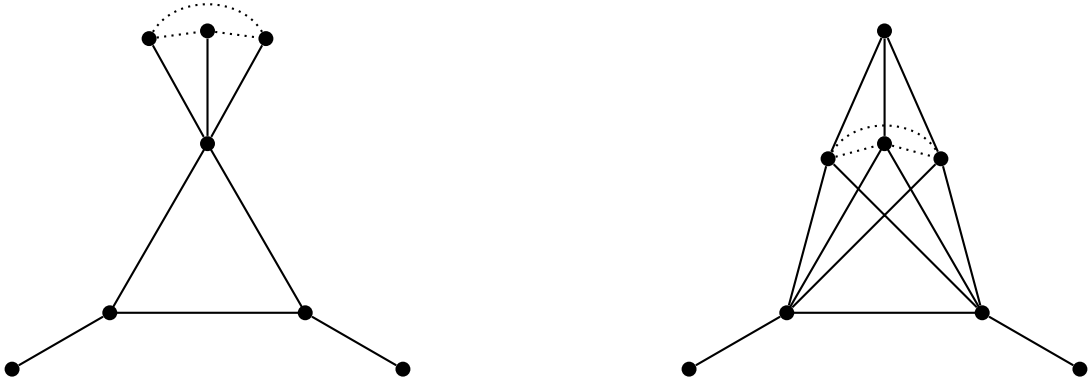
\begin{figure}[t] \centering
\begin{tikzpicture}[scale=1.5]
	\foreach \i in {0,1,2} {
		\node[circle,fill=black,scale=.5] (1i\i) at ({cos(90 + 120*\i)},{sin(90 + 120*\i)}) {};
		\node[circle,fill=black,scale=.5] (1o\i) at ({2*cos(90 + 120*\i)},{2*sin(90 + 120*\i)}) {};
		\node[circle,fill=black,scale=.5] (2i\i) at ({6+cos(90 + 120*\i)},{sin(90 + 120*\i)}) {};
		\node[circle,fill=black,scale=.5] (2o\i) at ({6+2*cos(90 + 120*\i)},{2*sin(90 + 120*\i)}) {};
		\draw[thick] (1i\i) -- (1o\i); \draw[thick] (2i\i) -- (2o\i);
		\draw[thick] ({cos(90 + 120*\i)},{sin(90 + 120*\i)}) -- ({cos(90 + 120*\i + 120)},{sin(90 + 120*\i + 120)});
		\draw[thick] ({6+cos(90 + 120*\i)},{sin(90 + 120*\i)}) -- ({6+cos(90 + 120*\i + 120)},{sin(90 + 120*\i + 120)});}
	\node[circle,fill=black,scale=.5] (c1) at ({2*cos(90+15)},{2*sin(90+15)}) {};
	\node[circle,fill=black,scale=.5] (c2) at ({2*cos(90-15)},{2*sin(90-15)}) {};
	\draw[thick] (c1) -- (1i0) (c2) -- (1i0); 
	\draw[thick,dotted] (c1) -- (1o0) -- (c2);
	\draw[thick,dotted] (c1) to[out=60,in=120] (c2);
	
	\node[circle,fill=black,scale=.5] (c1) at ({6+cos(90+30)},{sin(90+30)}) {};
	\node[circle,fill=black,scale=.5] (c2) at ({6+cos(90-30)},{sin(90-30)}) {};
	\draw[thick] (c1) -- (2o0) (c2) -- (2o0);
	\draw[thick,dotted] (c1) -- (2i0) -- (c2);
	\draw[thick] (2i1) -- (c1) -- (2i2) -- (c2) -- (2i1);
	\draw[thick,dotted] (c1) to[out=60,in=120] (c2);
\end{tikzpicture}
\caption{The two classes of $N_3$ graphs, with dotted lines to indicate potential edges.}\label{N22andN3}
\end{figure}

We show that extremal constructions are similar to nearly-balanced blow-ups of the net by investigating the density of two classes of 8-vertex graphs. The first is $N_{22}$, any 8-vertex graph obtained by cloning any two distinct vertices in the net, and the second is $N_3$, any 8-vertex graph obtained by cloning any one vertex twice in the net. For examples of $N_3$ graphs, see Figure \ref{N22andN3} where the dotted lines between vertices in the same part indicate that the edges are allowed but not necessarily present, as is standard procedure when drawing flags.

% Insert image of N22 and N3

Beginning our investigation of the extremal constructions by focusing on $N_{22}$ subgraphs will allow us to show (1) that much of the top-layer structure is present, and (2) the sizes of part sets are relatively balanced. We also investigate $N_3$ subgraphs to improve our bounds related to balanced part sizes at a small sacrifice to adherence to top-layer structure.

Now, we observe that in a balanced iterated blow-up of the net, $N_{22}(Z) - 5 \cdot N_{3}(Z)$ is computationally indistinguishable from 0 for any net $Z$, and so $N_{22}(Z) - a \cdot N_{3}(Z)$ should be a very small positive number when $5-a$ is a small positive number. (We will discuss the choice of $a$ later.) Further, in a balanced iterated blow-up of the net, no more than $1/36 + o(1)$ of the 8-vertex supergraphs of a net contained in a part are in either $N_{22}$ or $N_3$, so a net $Z$ which maximizes $N_{22}(Z) - a \cdot N_{3}(Z)$ should be a top-layer net. We use just one such net to classify the rest of the vertices in the graph into seven sets, six of which will function as approximations to the blow-up sets and one containing fewer than $0.2\%$ of the vertices which we cannot guarantee behave nicely.

We use flag algebra computations to bound $N_{22} - a \cdot N_3$, then apply an averaging argument to achieve a lower bound on the maximum of $N_{22}(Z) - a \cdot N_3(Z)$ over all nets $Z$. This bound is then included in quadratic programs designed to produce bounds on the number of vertices that we would need to delete, the number of edges that would need to be switched, and the relative sizes of the parts we find, and we solve these quadratic programs using Lagrange multipliers. Different choices for $a$ result in different sets of bounds, and we choose the set of bounds we need based on the arguments described next.

At this point, we will have argued that the extremal construction must be similar to a balanced blow-up of the net; some of the edges in the graph may disagree with the blow-up structure, the parts can vary in size to a small degree, and a very small set of vertices may simply not fit in any of the parts. To argue that no edges among the six primary parts disagree with the blow-up structure, we show that any pair of vertices violating the top-layer structure actually destroy more potential nets than they create. As such, we know that at least 99.8\% of the graph respects the top-layer structure. We then argue that the remaining 0.2\% of the vertices either cannot be placed in one of the six primary parts without creating conflict with at least 3\% of the vertices or else would unbalance the six parts too much. The bulk of the enumeration in the proof of Theorem \ref{net2} is found in these two arguments and is therefore a decent place to start when approaching the inducibility problem with a new graph.

We also know that all vertices must be in the same number of copies of the desired structure (Lemma \ref{AvgArg}), ignoring lower order terms. That means, in our problem, every vertex in an extremal construction is in $(24/1555 + o(1))\cdot\binom{n}{5} \approx 0.0001286 \cdot n^5$ nets. We then find an upper bound on the number of nets which contain one of the 0.2\% of vertices which create conflict with 3\% of the graph, and observe that this upper bound is less than $0.0001275 \cdot n^5$. In other words, all vertices in the graph are in $0.001286 \cdot n^5$ nets, but vertices which misbehave are in at most $0.0001275 \cdot n^5$ nets and hence do not exist in our extremal construction. This allows us to conclude that our extremal construction is a blow-up of the net graph, and we use two more arguments to show first that the parts are asymptotically balanced and finally that they differ by at most 1 vertex.

\section{Proof of Theorem 2}\label{ProofNotSketch}

\setcounter{thm}{-1}
\begin{prop}\label{prop3}
There exists some $a \in (4.9,5)$ and $n_0 \in \N$ such that every extremal graph $G$ on at least $n_0$ vertices satisfies: \begin{align*}
i(N,G) &< .017202164 \\
4 \cdot d(N_{22},G) - 3a \cdot d(N_3,G) &\geq 0.00071788399 + o(1) > 0.00071788398.
\end{align*}
\end{prop}

\noindent
\begin{proof}
This result follows from an application of the plain flag algebra method. We ran Flagmatic (\cite{flagmatic}) on 8 vertices which verified the first inequality for sufficiently large $n$. For the second inequality, we minimized the difference $4 d(N_{22},G) - 3a \cdot d(N_3,G)$, for several choices of $a \in (4.9,5)$, subject to the constraint that $i(N) > 24/1555$, the limiting density of the net in our conjectured extremal constructions. Our particular choice of $a$ which we carry forward into future calculations will be described in more detail later, but we note here that we will use $a = 4.99$. Values for $a$ closer to 5 and further from 5 gave collections of bounds that were less cooperative with future arguments, specifically Claim \ref{claim5}.
\end{proof}

Let $G$ be an extremal graph on $n$-vertices with $n$ large enough to satisfy the conditions of Proposition \ref{prop3}. Denote the set of all induced net in $G$ by $\mathcal{Z}$, let $a \in \R$ be taken from Proposition 3 (it will be specified later), and let $Z$ be an induced net maximizing $(N_{22}(Z) - a \cdot N_3(Z)) \binom{n-2}{4}$.

Below we obtain a lower bound on this linear combination in terms of flag algebra bounds by noting that it must be at least as large as the average value over all nets:
\begin{align*}
    \left(d(N_{22}(Z)) - a \cdot d(N_3(Z))\right)\binom{n - 6}{2} & \geq \frac{1}{|N|}\sum_{Y \in N}(d(N_{22}(Y)) - a d(N_3(Y))) \binom{n - 6}{2} \\
    &= \frac{4 d(N_{22}) - 3 a d(N_3) \binom{n}{8}}{d(N) \binom{n}{6}} \\
    &= \left(\frac{\frac{4}{28}d(N_{22}) - \frac{3a}{28}d(N_3)}{d(N)}\right)\binom{n-6}{2}.
\end{align*} Substituting the bounds from Proposition \ref{prop3} gives that
$$
d(N_{22}(Z)) - a \cdot d(N_3(Z)) > \frac{0.00071788398}{28 \cdot 0.17202164} = 0.00014904353852556826. % Label this inequality
$$
Label the vertices in $Z$ as $z_1,\dots,z_6$, and define sets of vertices $Z_i$ for $i \in [6]$ such that
$$
Z_i := \{v \in V(G) : G[(Z \backslash z_i) \cup v] \equiv N\}.
$$Effectively this creates sets of vertices $Z_i$ which look like $z_i$ to the other vertices in $Z$. For the sake of simplicity, let us label the vertices of $Z$ so that $E(Z) = \{z_iz_j: j = i+3 \text{ or } i,j\geq4\}$.

A funky pair $(x,y)$ is any pair of vertices $x \in Z_i$ and $y \in Z_j$ for $i \not= j$ where either $xy \in E(G)$ and $z_iz_j \not\in E(Z)$ or else $xy \not\in E(G)$ and $z_iz_j \in E(Z)$. In other words, a pair of vertices $x,y$ in distinct $Z_i,Z_j$ respectively is called funky if their adjacency status is inconsistent with that of $z_i$ and $z_j$. More specifically, we refer to the funky pair $x,y$ as a funky edge if $xy \in E(G)$ or as a funky non-edge if $xy \not\in E(G)$. We denote by $E_f$ the set of all funky pairs in $G$. It follows from the above inequality that
$$
\sum_{i=1}^3 |Z_i||Z_{i+3}| + \sum_{4 \leq i < j \leq 6} |Z_i||Z_j| - |E_f| - a \sum_{i\in[6]} |Z_i|^2/2
	> 0.000149043538 \cdot \binom{n-6}{2}.
$$
We then pick $X_i \subseteq Z_i$ and $X_0 = V(G) \backslash \op{\ds\cup_{i=1}^6 X_i}$ so as to maximize the left-hand side of (4.1), below, where we also define $x_i := \frac{1}{n}|X_i|$ for $i = 0,\dots,6$ and $f$ as the density of funky pairs with neither end in $X_0$, normalized by a factor $\binom{n}{2}^{-1}$. This gives the following useful bound to be used as a constraint in the proof of Claim \ref{claim4}:
\begin{align}\label{inequality3}
2 \sum_{1 \leq i < j \leq 6} x_ix_j - 2f - a\sum_{i=1}^6 x_i^2 > 0.000149043538.
\end{align}

In the following claim, we specify a specific real number $a$, which gives corresponding bounds. Note that different choices for $a$ will result in different inequalities below and that the choice of $a$ is primarily made in such a way to allow arguments \emph{later in the proof} to go through. Specifically, we chose $a = 4.99$ so that arguments in Claim \ref{claim5} produce the desired contradictions.

\begin{claim}\label{claim4}
For $a = 4.99$, we have for all $i \in [6]$ that
\begin{align*}
0.165791592261 \leq x_i &\leq 0.167541741072, \\
x_0 &\leq 0.00165262197319, \\
f &\leq 0.0000027521.
\end{align*}
\end{claim}

\noindent
\begin{proof}
These bounds are achieved via a straight-forward application of quadratic programming and Lagrange multipliers; see Appendix A for details. The structure of this proof follows closely to that of Claim 4 in \cite{balHuLidPfen}.
\end{proof}

Define $\xmin := 0.165791592261$ and $\xmax := 0.167541741072$. We can show that the funky degree of any $x \in V(G) \backslash X_0$ satisfies $d_f(x) \leq 1 - (1+a) \cdot x' \approx 0.0069084$: if we move $v$ from $X_1$ to $X_0$, then the left hand side of \ref{inequality3} will decrease by
$$
\frac{1}{n}(2(x_2 + \cdots + x_6) - 2d_f(v) - 2 \cdot a \cdot x_1 + o(1)).
$$This quantity must be positive as $X_0, X_1, \dots, X_6$ were chosen to maximize the left hand side of \ref{inequality3}. This together with the bounds on $x_i$ from Claim \ref{claim4} implies that
$$
d_f(v) \leq x_2 + \cdots + x_6 - d_f(v) - a \cdot x_1 + o(1) \leq 1 - (1+a) \cdot x' + o(1).
$$

Given that the funky degree of vertices in $V(G) \backslash X_0$ is bounded above by $0.0069084$, we now show that all funky pairs in $G$ must involve a vertex in $X_0$. This is the first place in which our strategy diverges from that of Balogh, Hu, Lidick\'{y}, and Pfender in \cite{balHuLidPfen}. The $C_5$ problem involves a great deal of symmetry as $C_5$ is both vertex-transitive and self-complementary. The net graph has neither of these properties, although we will certainly take advantage of the symmetries that the net graph \emph{does} have. Instead, we will iteratively show that certain types of funky pairs which avoid $X_0$ do not exist, and then use these observations to reduce the number of possible nets in which the more troublesome funky pairs may reside.

In the proof of the following claim, we will use some unorthodox terminology, which we now introduce. We refer to $X_1,\dots,X_6$ as \emph{blobs}, and we specify $X_1$, $X_2$, and $X_3$ as \emph{outer blobs} as a reference to the vertices $z_1$, $z_2$, and $z_3$ in $Z$, and we specify $X_4$, $X_5$, and $X_6$ as \emph{inner blobs} for an equivalent reason. We define the \emph{blob distance} of two vertices $x \in X_i$ and $y \in X_j$ for $i,j \in [6]$ as the distance of $z_i$ and $z_j$ in $Z$, we say that two or more vertices are \emph{coblobular} if their blob distance is 0. We say that a subgraph of $G$ is \emph{blob-induced} or \emph{respects the blob-structure} if the subgraph contains no funky edges, usually with no pair of coblobular vertices.

\begin{claim}\label{claim5}
There are no funky pairs in $X_1 \cup \cdots \cup X_6$.
\end{claim}

\noindent
\begin{proof}
We prove this claim by contradiction. Let $u,v \in (X_1 \cup \cdots \cup X_6)$ be a funky pair in $G$, and let $G'$ be obtained from $G$ by changing the status of the edge $uv$. We will show that the bounds in Claim 4 imply that $G'$ must contain more nets than the extremal construction $G$.

Using the bounds in Claim \ref{claim4}, we can find at least $0.00069 \cdot n^4$ nets containing the pair $u,v$ in $G'.$ Without loss of generality, the funky pair $u,v$ are of one of the following four types: (1) $u \in X_4$, $v \in X_5$, (2) $u \in X_1$, $v \in X_4$, (3) $u \in X_1$, $v \in X_5$, or (4) $u \in X_1$, $v \in X_2$. We will show that if $u,v$ are a funky pair of one of these types, in this order, then $G$ has fewer than $0.00069 \cdot n^4$ nets which contain the pair $u,v$.

First, we note that there are at most $x_0/6 \cdot n^4$ nets using at least one vertex in $X_0$. Second, we count the number of nets in $G$ containing $u,v$ and at least two other vertices $u',v'$ in funky pairs by observing that there are at most $f/2 \cdot n^4$ nets wherein $(u',v')$ is funky and at most $(d_f(u) + d_f(v))/2 \cdot n^4$ where $u'$ and $v'$ are each funky to at least one of $u$ or $v$.

Third, we want to bound the number of nets wherein $(u,v)$ is the only funky pair. Let $N$ be a net in $G$ with $u,v \in V(N)$ and note that $N$ contains a $P_4$ which itself has no funky pairs. As $P_4$ is prime, we note that this $P_4$ is either blob-induced or coblobular.

Suppose towards contradiction that this $P_4$ is blob-induced and assume without loss of generality that it lies in $X_1$, $X_4$, $X_6$, and $X_3$. If $u,v$ are not in this $P_4$, then they must be blob-distance at least 2 from each other, forcing at least one of them to have at least one additional funky partner within the $P_4$, a contradiction. Otherwise, there exists a second non-funky $P_4$ in $N$ which must then be blob-induced; any placement of this non-funky $P_4$ contradicts the assumption that $(u,v)$ is the only funky pair in $N$.

Therefore, the original $P_4$ is coblobular, which implies first that the third triangle vertex must be coblobular with the $P_4$. Since the net has \emph{some} funky pair, it follows that $u,v$ must be a funky pendant-edge in $N$; as such, we have at most $(\xmax n)^4/12$ nets in $G$ where $u,v$ is the only funky pair.

The nets in $G$ that have yet to be counted are those containing precisely one other vertex $w$ which is funky to $u$, $v$, or both; there are at most $c \cdot (d_fn) \cdot (\xmax n)^3$ for some constant $c$. In order to arrive at our desired contradiction, we would need $c$ to satisfy
$$
(\xmax)^4/12 + c \cdot d_f \cdot (\xmax)^3 + x_0/6 + d^2_f + f/2 < 0.00069,
$$which is the case when $c < 9.522$. \\

\noi
\emph{Sub-claim 5.1}
There are no funky pairs in $X_4 \cup X_5 \cup X_6$.

\noindent
\begin{proof}
Let us first enumerate the nets containing $u$, $v$, and a vertex $w$ which is funky with both $u$ and $v$, and suppose without loss of generality that $u \in X_4$ and $v \in X_5$.

There are no nets in which $w$ is in $X_1 \cup X_2$ by the following contradiction argument. Suppose $N$ is just such a net and $w \in X_1$ without loss of generality; this implies that a mutual non-neighbor $x \in V(N)$ of $v$ and $w$ must be located in $X_3$, and the existence of a $v-x$ path in $N$ implies the existence of $y \in V(N) \cap X_6$. Clearly $y$ is a triangle vertex in $N$ but its neighborhood in $N$ forms an independent set, a contradiction.

There are no nets in which $w$ is in $X_6$ by the following contradiction argument. Suppose $N$ is just such a net and observe that $X_6 \cap V(N)$ is empty as $N$ contains no $C_4$s. We can then see that $w$ is a triangle vertex, one of $u$ or $v$ is a triangle vertex, but the third triangle vertex cannot be adjacent to the other two without additional funky pairs.

Now we count the nets with $w \in X_6$. We know that at most one of $u$, $v$, $w$ is a triangle vertex, and so the remaining two vertices $t_1,t_2$ must be blob-distance at most one from each of $u$, $v$, and $w$; hence $t_1$ and $t_2$ are in different interior blobs as nets do not contain $C_4$s. The last vertex is then a pendant of the triangle vertex in $\{u,v,w\}$, and so must be blob-distance at least 2 from both $t_1$ and $t_2$. There are 3 distinct configurations for $t_1$ and $t_2$, which specifies precisely the location of the last pendant, implying that there are at most $3d_f(\xmax)^3n^4$ such nets.

Lastly, we count the nets wherein $w$ is funky only with $u$ and apply symmetry to acquire the number of nets wherein $w$ is funky only with $v$. As $uv \not\in E(G)$, we have that within any such net, $v$ is contained in a non-funky $P_4$. If this $P_4$ is completely contained in $X_5$, then $u$ is a triangle vertex, $v$ is a pendant, and $w$ is also a pendant in this non-funky $P_4$. The pendant of $u$ can either be in $X_1$ or $X_5$, implying that there are $\op{\frac{1}{2} + \frac{1}{6}}d_f\xmax^3n^4$ such nets. Otherwise, $v$ is a triangle vertex, $p_v \in X_2$, and $u$ is a pendant vertex. Either the $P_4$ intersects $X_4$ or $X_6$, the latter implying that $w \in X_6$, $p_2 \in X_3$, and $t_u \in X_4$ (giving $d_f\xmax^3n^4$ nets), and the former implying that $w \in X_1$, $t_w \in X_4$, and $t_u \in X_6$ (giving $d_f\xmax^3n^4$ nets). So, we have $\frac{8}{3}d_f\xmax^3n^4$ nets wherein $w$ is funky only with $u$, and thus $(3 + \frac{16}{3})d_f\xmax^3n^4 < 9.522d_f\xmax^3n^4$. This gives the desired contradiction, proving Sub-claim 5.1.
\end{proof}

\noi
\emph{Sub-claim 5.2}
There are no funky non-edges in $X_1 \cup \cdots \cup X_6$.

\noindent
\begin{proof}
First note that Sub-claim 5.1 implies that the only possible funky non-edges exist between an outer blob and its adjacent inner blob, so we will suppose towards contradiction and without loss of generality that $u \in X_1$ and $v \in X_4$.

Let us first determine how many nets containing $u$, $v$, and $w$ wherein $w$ is funky to both $u$ and $v$. By Sub-claim 1, we know that $w \not\in X_5 \cup X_6$. Therefore, $w \in X_2 \cup X_3$. No vertex in any such net can be adjacent to both $u$ and $w$, so $u = p_w$ and $v$ is a triangle-vertex. If $w \in X_2$, then we need a third triangle vertex $x$ in $X_5$, and then $p_x \in X_3$. A symmetric argument works for $w \in X_3$, so we have $2d_f(x'')^3n^4$ nets of this form.

Now let us count the number of nets containing $u$, $v$, and $w$ wherein $w$ is funky only to $v$, and note by Sub-claim 1 that none exist where $w \in X_5 \cup X_6.$ Furthermore, we know that any such net contains a non-funky $P_4$ which excludes $v$ but includes $u$. If this $P_4$ is contained in $O1$, then $w \in X_1$, both $u$ and $w$ are pendants, $v$ is a triangle vertex, and its pendant can be anywhere in $X_1 \cup X_5 \cup X_6$, giving a count of $(\frac{1}{6} + 2 \cdot \frac{1}{2})d_f\xmax^3n^4$ nets. Otherwise, the $P_4$ is blob-induced and $u$ is a pendant of a vertex $t_u \not = v$ in $X_4$. By symmetry, we assume that the $P_4$ intersects $X_2$, and note that $v$ must be adjacent to the vertex $y \in X_5$ from the $P_4$. Since $p_y \in X_2$, it follows that $v$ is a triangle vertex and $w = p_v \in X_3$, giving a count of $d_f\xmax^3n^4$ nets. Hence we have $\frac{19}{6}d_f\xmax^3n^4$ nets wherein $w$ is funky only $v$.

Finally, we count the number of nets containing $u$, $v$, and $w$ wherein $w$ is funky only to $u$. Observe that any such net contains a non-funky $P_4$ containing $v$. If this $P_4$ is entirely within $X_4$, then $u$ is a triangle vertex, $v$ and $w$ are the pendants in this $P_4$, and $p_u$ must also be contained in $X_4$ due to adjacency restrictions, giving a count of $\frac{1}{6}d_f\xmax^3n^4$ nets. Otherwise, the $P_4$ is blob-induced, $v$ is a triangle vertex, $p_v \not= u$ is contained in $X_1$, $u$ is a pendant, the triangle is blob-induced, and $w = t_u$ is not the triangle vertex in the aforementioned $P_4$, giving a count of $2d_f\xmax^3n^4$ nets by symmetry. This results in $\frac{13}{6}d_f\xmax^3n^4$ nets wherein $w$ is funky only to $u$.

Therefore, there are at most $(\frac{13}{6} + \frac{19}{6} + 2)d_f\xmax^3n^4$ nets containing the funky pair $u \in X_1$ and $v \in X_4$. Since $\frac{13}{6} + \frac{19}{6} + 2 < 9.522$, we have our contradiction. This implies not only that no funky pair exists in $X_1 \cup X_4$, but further that no funky non-edge exists in $X_1 \cup \cdots \cup X_6$, proving Sub-claim 5.2.
\end{proof}

\noi
\emph{Sub-claim 5.3}
There are no funky pairs in $X_1 \cup X_5$.

\noindent
\begin{proof}
Without loss of generality, suppose $u \in X_1$ and $v \in X_5$. We will begin by counting the nets wherein $w$ is funky to both $u$ and $v$, noting that $w$ is therefore contained in $X_3$. So we know that the triangle in any such net is necessarily formed by $u$, $v$, and $w$; the pendants of both $u$ and $w$ must be in $X_1$ and $X_3$, respectively, while $p_v$ may be anywhere in $X_5 \cup X_2$, giving $2d_f\xmax^3n^4$ nets.

We now count the number of nets wherein $w$ is funky only to $u$. By Sub-claim 5.2, we have that $w \not\in X_4$. We get $d_f\xmax^3n^4$ nets when $w \in X_6$ and no nets when $w \in X_2$. When $w \in X_5$, the only possible location for a third neighbor of $u$ is in $X_1$, implying both $v$ and $w$ are triangle vertices; their pendants cannot be blob-adjacent, so $p_v,p_w \in X_5$ giving $\frac{1}{2}d_f\xmax^3n^4$ nets. The only remaining case is when $w \in X_3:$ as nets do not contain $C_4$s, $X_6$ is empty, $w$ would have to be the pendant of $u$, we would need a triangle vertex $t \in X_4$, $p_t \in X_1$, and $p_v \in X_2$, giving a count of $d_f\xmax^3n^4$ nets. This implies that there are at most $\frac{5}{2}d_f\xmax^3n^4$ nets wherein $w$ is funky only to $u$.

Finally, we count the number of nets wherein $w$ is funky only to $v$. Sub-claim 5.2 implies that $w \in X_1 \cup X_3$, so we first count those where $w \in X_1$. In this case, no nets intersect $X_4$ since nets do not contain $C_4$s, and hence the non-funky $P_4$ which does not contain $v$ must be entirely within $X_1$. As $v$ cannot have three neighbors in $X_1$, we have that $p_v$ is anywhere in $X_2 \cup X_5 \cup X_6$, giving a count of $3\cdot\frac{1}{2}d_f\xmax^3n^4$ nets. Next, we see that if $w \in X_3$, then either $u$ or $w$ is a triangle vertex, but not both. In either case, their pendant would have to share their blob, which then contradicts the existence of a third triangle vertex. So we have at most $\frac{3}{2}d_f\xmax^3n^4$ nets wherein $w$ is funky only to $u$.

As $2 + \frac{5}{2} + \frac{3}{2} < 9.522$, we have our desired contradiction, proving Sub-claim 5.3.
\end{proof}

Taking the above three sub-claims together, we have that all funky pairs in $X1 \cup \cdots \cup X_6$ must be between outer blobs. To finish the proof of Claim \ref{claim5}, we suppose toward the familiar contradiction, without loss of generality, that $u \in X_1$ and $v \in X_2$. We will count the number of nets containing $u$, $v$, and $w$, wherein $w$ is funky to at least one of $u$ and $v$.

We begin by counting the nets wherein $w$ is funky to both $u$ and $v$, noting that our sub-claims imply that $w \in X_3$. In this case, there can be at most one vertex in an interior blob, giving a count of $4d_f\xmax^3n^4$ nets.

We will then use symmetry to only count the nets wherein $w$ is funky only to $u$. We know that the third neighbor $x$ of $u$ in any such net is contained in $X_1 \cup X_4$, implying that $w$ must be a neighbor of $v$. Therefore, we have that $w \in X_2$ with $v$, $p_v$ and $p_w$, giving at most $2d_f\xmax^3n^4$ nets in this case.

Doubling our most recent case to account for the nets wherein $w$ is funky only to $v$, we have that there are at most $8d_f\xmax^3n^3$ nets containing $u$, $v$, and some $w$ which is funky to one or both of $u$ and $v$, giving our desired contradiction and completing the proof that there are no funky pairs in $X_1 \cup \cdots \cup X_6$.
\end{proof}

We now intend to show that in an extremal example, the trash is empty (i.e. $x_0 = 0$). We will do this by considering the ramifications of moving a vertex $x \in X_0$ to a blob, specifically as they relate to $d_f(x)$. Claim \ref{claim6} gives a lower bound on the potential funky degree of trash vertices, while Claim \ref{claim7} shows that every vertex in an extremal example must be in roughly the same number of nets. Claim \ref{claim8} then argues that trash vertices, by virtue of having high funky degree, cannot be in the correct number of nets. This will imply that $X_0$ contains only vertices which can be placed in a blob without creating any funky pairs.

\begin{claim}\label{claim6}
If a vertex $x \in X_0$ is moved to $X_i$, then either $d_f(x) = 0$ or else
$$
d_f(x) \geq \left\{ \begin{matrix}
	0.0433316 & i \in [3] \\
	0.0322447 & i \not\in [3].
\end{matrix} \right.
$$
\end{claim}

\noindent
\begin{proof}
Prior to moving $x$ into a blob, there were no funky pairs; as such, after moving $x$, all funky pairs must include $x$. Suppose that $d_f(x) > 0$ and let $w$ be a funky partner of $x$. For the purposes of case analysis, there are six orientations for the funky pair $(x,w)$. For the sake of brevity, we will only show the arguments for one of the cases here, and the others can be found in the Appendix B.

To achieve a lower-bound on $d_f(x)$, we will take an approach similar to the proof of Claim \ref{claim5}. As all funky pairs include $x$, we know that every net containing the funky pair $(x,w)$ contains at least one non-funky $P_4$. Further, we need only find the placement of the $d_f(x)n$ funky partners of $x$ which maximizes the number of nets containing $(x,w)$. Regardless of the location of $x$ and $w$, we can get a trivial upper bound on the number of nets containing both $x$ and another vertex in the trash as $\frac{n^3}{6}(x_0n)$.

The case we will investigate is when $x \in X_1$ and $w \in X_2$, and we begin by determining all possible configurations of nets containing the edge $xw$. First, let us determine those configurations wherein all of the other 4 vertices are coblobular: \begin{itemize}
	\item[(a)] a non-funky $P_4$ is contained in $X_1$,
	\item[(b)] four vertices lie in $X_2$ and $x = p_w$,
	\item[(c)] four vertices lie in $X_2$ and $w = p_x$,
	\item[(d)] four vertices lie in $X_2$ and $wx \in \Delta$, or
	\item[(e)] four vertices lie in one of $X_3$, $X_4$ or $X_6$.
\end{itemize}Note that four vertices cannot lie in $X_5$ as $w$ can neither be funky to any other vertex in the graph nor have degree 5. Together, constructions (a) and (b) contribute at most $2(\xmax n)^4/24$ nets. Notice, however, that constructions (c)-(e) all require the choice of 2 additional funky partners of $x$; we could take a trivial upper bound of $5\cdot\frac{(d_f(x)n)^2}{2}\cdot\frac{(\xmax n)^2}{2}$ but the entire funky neighborhood of $x$ cannot be in 4 places at once. We will leverage this observation shortly.

If there is a non-funky blob-induced $P_4$ in the net, then it must neither contain $w$ nor intersect $X_5$. Therefore, any such net must have its non-funky blob-induced $P_4$ in $X_1$, $X_4$, $X_6$, and $X_3$. In this case, $x = t_w$ and contributes at most $(\xmax n)^3$ for each funky neighbor in $X_6$. The only remaining case is that the non-funky $P_4$ contains $w$ and lies in $X_2$, implying that $xw \in \Delta$ and $p_x$ is in $X_1 \cup X_4 \cup X_6 \cup X_3$.

So, there are at most $(\xmax n)^4/12$ nets containing $xw$ and no other funky partners of $x$. We now take advantage of the observation that the funky neighbors of $x$ cannot be in more than one place at a time. Define $d_i(x)$ to be the density of funky partners in $X_i$. First, we determine an upper bound on the nets containing $xw$ and one other funky partner of $x$: \begin{align*}
(d_6(x)n)\cdot(\xmax n)^3 + (d_2(x)n) \cdot \frac{(\xmax n)^2}{2} \cdot (\xmax n)
	&= \xmax^3n^4\op{d_6(x) + \frac{1}{2}d_2(x)} \\
	& \leq d_f(x)\xmax^3n^4,
\end{align*}which is achieved if all funky partners of $x$ lie in $X_6$. Next, we determine an upper bound on the nets containing $xw$ and two other funky partners of $x$: \begin{align*}
n^4\frac{\xmax^2}{2}\op{\frac{(d_2(x)^2 + d_3(x)^2 + d_4(x)^2 + d_6(x)^2)}{2} + d_2(x)(d_3(x) + d_6(x))}.
\end{align*}We observe that any maximum of this expression will occur when $d_4(x) = 0$. Further, if $d_i(x) > 0$ for each $i \in \{2,3,6\}$, then there exist pairs of funky partners $(v_3,v_6) \in X_3 \times X_6$ which do not contribute to the net count; as such, a maximum of our expression occurs when $d_6(x) = 0$ as well. Therefore, we have that the number of nets containing $xw$ and two other funky partners of $x$ is bounded above by
\begin{align*}
n^4\frac{\xmax^2}{2}\op{\frac{(d_2(x)^2 + d_3(x)^2)}{2} + d_2(x)d_3(x)}
 &= n^4\frac{\xmax^2}{2} \cdot \frac{1}{2}(d_2(x) + d_3(x))^2 \\
 &\leq n^4\frac{\xmax^2}{2} \cdot \frac{1}{2}(d_f(x))^2 \\
 &= \frac{1}{4}d_f(x)^2\xmax^2n^4.
\end{align*} Therefore, it follows that when $x \in X_1$ and $w \in X_2$, the number of nets containing $xw$ is bounded above by
$$
n^4 \op{\frac{\xmax^4}{12} + d_f(x)\xmax^3 + \frac{1}{4}d_f(x)^2\xmax^2 + \frac{x_0}{6}}.
$$If we delete the edge $xw$, we have at least $\xmin^3(\xmin - d_f(x))n^4$ nets. As $G$ is extremal, it follows that
$$
\xmin^3(\xmin - d_f(x)) \leq \frac{\xmax^4}{12} + d_f(x)\xmax^3 + \frac{1}{4}d_f(x)^2\xmax^2 + \frac{x_0}{6},
$$which implies that $d_f(x) \geq 0.0433316$. Following similar arguments, detailed in Appendix B, we arrive at the bounds found in Table \ref{netBounds}, concluding the proof of Claim \ref{claim6}.

\begin{table}[t]
\centering
\caption{Number of copies of nets containing $x,w$ in given config class.}\label{netBounds}
\begin{tabular}{|c|c|c|c|}
\hline
$x$ & $w$ & Upper Bounds for Nets & Lower Bound for $d_f(x)$ \\
\hline
$X_1$ & $X_2$ & $\frac{1}{6}x_0 + \frac{1}{12}\xmax^4 + d_f(x)\xmax^3 + \frac{1}{4}d_f(x)^2\xmax^2$ & 0.0433316 \\ & & & \\
$X_1$ & $X_4$ & $\frac{1}{6}x_0 + \frac{1}{2}d_f(x)\xmax^3 + \frac{1}{2}d_f(x)^2\xmax^2 + \frac{1}{6}d_f(x)^3\xmax$ & 0.0610118 \\ & & & \\
$X_1$ & $X_5$ & $\frac{1}{6}x_0 + \frac{1}{12}\xmax^4 + d_f(x)\xmax^3 + \frac{1}{4}d_f(x)^2\xmax^2$ & 0.0433316 \\ & & & \\
$X_4$ & $X_1$ & $\frac{1}{6}x_0 + \frac{13}{6}d_f(x)\xmax^3 + \frac{1}{8}d_f(x)^2\xmax^2 + \frac{1}{6}d_f(x)^3\xmax$ & 0.0322447 \\ & & & \\
$X_4$ & $X_2$ & $\frac{1}{6}x_0 + \frac{1}{12}\xmax^4 + \frac{3}{2}d_f(x)\xmax^3 + \frac{1}{4}d_f(x)^2\xmax^2$ & 0.0349529 \\ & & & \\
$X_4$ & $X_5$ & $\frac{1}{6}x_0 + d_f(x)\xmax^3 + \frac{1}{8}d_f(x)^2\xmax^2 + \frac{1}{6}d_f(x)^3\xmax$ & 0.0504913\\ \hline
\end{tabular}
\end{table}
\end{proof}

The following claim follows directly from Lemma \ref{AvgArg} and is placed here primarily to remind the reader ahead of Claim \ref{claim8} in roughly how many nets any given vertex in $G$ must reside.

\begin{claim}\label{claim7}
Every vertex of the extremal graph $G$ is contained in at least $(24/1555 + o(1))\binom{n}{5} \approx 0.00012861736n^5$ induced net graphs.
\end{claim}

The following claim will then imply that our extremal example $G$ is a (not necessarily balanced or iterated) blow-up of the net graph.

\begin{claim}\label{claim8}
For each $v \in X_0$, there is some $i \in [6]$ so that $v$ can be moved to $X_i$ without creating any funky pairs.
\end{claim}

\noindent
\begin{proof}
Assume there exists some $x \in X_0$ so that $d_f(x) > 0$ if $x$ were to be placed in any of the 6 main blobs. We will show that $x$ is in far fewer than $(24/1555 + o(1))\binom{n}{5}$ nets, a contradiction.

Let $a_in$ be the number of neighbors of $x$ in $X_i$ and $b_in$ be the number of neighbors of $x$ in $X_i$ for each $i \in \{0,1,\dots,6\}$. Normalizing by $n^4$, we will let $A$ be the number of nets containing $x$ and vertices from 5 different blobs, $B$ be the number of nets containing $x$ and 5 coblobular vertices, and $C$ be all other nets not involving trash vertices.

We have that
\begin{align*}
A &\leq b_2b_3a_4b_5b_6 + b_1b_3b_4a_5b_6 + b_1b_2b_4b_5a_6 
	+ a_1b_2b_3a_5a_6 + b_1a_2b_3a_4a_6 + b_1b_2a_3a_4a_5, \\
B &\leq \sum_{i=1}^6\op{\frac{a_i^3b_i^2}{3!2!} + \frac{a_ib_i^4}{1!4!}}\text{, and } \\
C &\leq \frac{1}{2!2!}\sum_{i=1}^3a_i^2b_i^2\op{-a_i-a_{i+3} + \sum_{1 \leq j \leq 6}a_j}
	+ \frac{1}{2!2!}\sum_{i=4}^6a_i^2b_i^2(a_1+a_2+a_3-a_{i-3}).
\end{align*}A convenient way to convince yourself that these counts are correct is to observe that every net must contain either a non-funky bull or a non-funky $P_4 + K_1$, and note that these structures are either blob-induced or component-wise coblobular and respect the blob-structure. We will say that $D$ is the number of nets containing $x$ and at least one other trash vertex.

The goal of this claim, then, is to show that the following program is bounded away from $(24/1555)/5! = 1/7775 \approx 0.000128617363$:
$$
(P) \left\{ \begin{matrix*}[l]
\text{maximize} & A + B + C + D \\
\text{subject to} & \sum_{i=0}^6(a_i+b_i) = 1 \\
	& (x') \leq a_i + b_i \leq (x'') \text{ for } i \in [6], \\
	& a_0 + b_0 \leq x_0, \\
	& a_2 + a_3 + b_4 + a_5 + a_6 \geq 0.0433316, \\
	& a_1 + a_3 + a_4 + b_5 + a_5 \geq 0.0433316, \\
	& a_1 + a_2 + a_4 + a_5 + b_6 \geq 0.0433316, \\
	& b_1 + a_2 + a_3 + b_5 + b_6 \geq 0.0322447, \\
	& a_1 + b_2 + a_3 + b_4 + b_6 \geq 0.0322447, \\
	& a_1 + a_2 + b_3 + b_4 + b_5 \geq 0.0322447, \\
	& a_i,b_i \geq 0 \text{ for } i \in \{0,1,\dots,6\}.
\end{matrix*} \right.
$$

However, we do not have a form for $D$ as we will be optimizing over all possible ways to place trash vertices into the blobs and assume that the trash vertices always contribute to nets a vertex from its blob contributes to nets. So we relax the above problem the following:
$$
(P') \left\{ \begin{matrix*}[l]
\text{maximize} & f=A + B + C \\
\text{subject to} & \sum_{i=1}^6(a_i+b_i) = 1 \\
	& (x') \leq a_i + b_i \leq (x'' + x_0) \text{ for } i \in [6], \\
	& a_2 + a_3 + b_4 + a_5 + a_6 \geq 0.0433316, \\
	& a_1 + a_3 + a_4 + b_5 + a_5 \geq 0.0433316, \\
	& a_1 + a_2 + a_4 + a_5 + b_6 \geq 0.0433316, \\
	& b_1 + a_2 + a_3 + b_5 + b_6 \geq 0.0322447, \\
	& a_1 + b_2 + a_3 + b_4 + b_6 \geq 0.0322447, \\
	& a_1 + a_2 + b_3 + b_4 + b_5 \geq 0.0322447, \\
	& a_i,b_i \geq 0 \text{ for } i \in \{0,1,\dots,6\}.
\end{matrix*} \right.
$$
We will discretize the space of possible solutions to $(P')$, determine the value of the objective function at the center of each cell, and use a bound on the gradient to show that the function is bounded above by $0.0001275$ in each cell. If the global bound on the gradient is not sufficient in bounding the optimization function, we generate a bound on the gradient within the cell, and if necessary, refine the discretization within the cell. For every $a_i$ and $b_i$, we check $s+1 = 51$ equally spaced values between 0 and $x'' + x_0$, inclusive. By this, we have a grid of $s^{12}$ boxes where every feasible solution of $(P')$, and hence of $(P)$, is in one of the boxes.

To determine a global bound on the gradient, we find the partial derivatives of $f$:
\begin{align*}
\frac{\partial f}{\partial a_1}
	&= b_2b_3a_5a_6 + \frac{1}{24}b_1^4 + \frac{1}{4}a_1^2b_1^2 + \frac{1}{2}(a_1b_1^2)(a_2+a_3+a_5+a_6) + \frac{1}{4}(a_2^2b_2^2 + a_3^2b_3^2 + a_5^2b_5^2 + a_6^2b_6^2) \\
\frac{\partial f}{\partial a_4}
	&= b_2b_3b_5b_6 + b_1a_2b_3a_6 + b_1b_2a_3a_5 + \frac{1}{24}b_4^4 + \frac{1}{4}a_4^2b_4^2	
		+ \frac{1}{4}(a_2^2b_2^2 + a_3^2b_3^2) + \frac{1}{2}(a_4b_4^2)(a_2+a_3) \\
\frac{\partial f}{\partial b_1}
	&= b_3b_4a_5b_6 + b_2b_4b_5a_6 + a_2b_3a_4a_6 + b_2a_3a_4a_5 + \frac{1}{6}a_1^3b_1 + \frac{1}{6}a_1b_1^3 + \frac{1}{2}a_1^2b_1(a_2 + a_3 + a_5 + a_6) \\
\frac{\partial f}{\partial b_4}
	&= b_1b_3a_5b_6 + b_1b_2b_5a_6 + \frac{1}{6}a_4^3b_4 + \frac{1}{6}a_4b_4^3 + \frac{1}{2}a_4^2b_4(a_2+a_3).
\end{align*}We were able to bound each of these partials by $\frac{4}{3}(x'' + x_0)^4$. To do so, one should observe the portion of the partial contributed by $A$ can be bounded above by $(x'' + x_0)^4$ in each case, and therefore the rest of the partial can be bounded by $(1/3)(x'' + x_0)^3$ by taking advantage of obvious symmetries. One should note that the partial taken with respect to $b_1$ is the closest to meeting our bound, while the others are closer to $\frac{7}{6}(x'' + x_0)^3$.

As each 12-dimensional cell has side-length $1/s$, the objective function can exceed the value at its center by at most $12 \cdot \frac{1/2}{2} \cdot \frac{4}{3}(x'' + x_0)^4$. The local bounds on the gradient are obtained in our algorithm are achieved by substituting $(a_i + (1/t)/2)$ and $(b_i + (1/t)/2)$ into each partial, where $t$ is the side-length of the current cell which may or may not be refined, and we simply take the steepest direction of ascent as our bound to replace $\frac{4}{3}(x'' + x_0)^4$.

Using $s = 50$, we successfully bounded the objective function below $0.0001275$; a handful of cells required refinement, but not one cell required more than 6 refinements. The C++ code can be found at \url{https://sites.google.com/view/ablumenthal/}.

As the number of nets containing $x$ is bounded away from the average, we have contradicted the existence of a trash vertex which cannot be placed in any blob without creating funky pairs, completing the proof of Claim 8.
\end{proof}

Since every trash vertex can be placed into one of the six blobs without creating any funky pairs with vertices originally in $X_1 \cup \cdots \cup X_6$. Therefore, we simply add each trash vertex into its corresponding blob; at worst, we may have funky pairs involving trash vertices in different blobs, but we simply apply Claim \ref{claim5} noting that our bounds on $d_f(x)$ and $f$ are even more strict.

Therefore, we have established the top-layer structure of $G$ where each part is roughly one-sixth of the graph. As such, an induced net in $G$ can only be found by picking a single vertex from each of the 6 blobs or by picking all 6 within the same blob, implying that
$$
N(n) = (x_1 \cdots x_6)n^6 + N(x_1n) + \cdots + N(x_6n).
$$By averaging over all subgraphs of $G$ of order $n-1$, we have that $N(n)\leq \frac{n}{n-6}N(n-1)$ for all $n$, so
$$
\ell := \lim_{n\to\infty}\frac{N(n)}{\binom{n}{6}}
$$exists and so satisfies
$$
\ell + o(1) = 6!(x_1 \cdots x_6) + \ell \cdot (x_1^6 + \cdots x_6^6),
$$which implies that that $x_i = \frac{1}{6} + o(1)$ and $\ell = \frac{24}{1555}$, giving the constraints on the $x_i$'s. In order to complete the proof of Theorem \ref{net2}, we need only show that $|X_i| - |X_j| \leq 1$ for all $i,j \in [6]$.

\begin{claim}\label{claim9}
For $n$ large enough, we have $|X_i| - |X_j| \leq 1$ for all $i,j \in [6]$.
\end{claim}

\noindent
\begin{proof}
Let $i,j \in [6]$ satisfy $|X_i| - |X_j| \geq 2$, noting that $i \not= j$. Let $v \in X_i$ be chosen to minimize $N^u$ over vertices in $X_i$ and $w \in X_j$ be chosen to maximize $N^w$ over vertices in $X_j$. As $G$ is extremal, $N^v + N^{v,w} - N^w \geq 0$; otherwise, we can increase the number of nets by replacing $v$ with a copy of $w$.

Let $y_i :=|X_i| = x_in$ and $y_j := |X_j| = x_jn$. By the monotonicity of $\frac{N(n)}{n^6}$, we have
$$
\frac{24}{1555} + o(1) \geq \frac{N(y_j)}{\binom{y_j}{6}} \geq \frac{N(y_i)}{\binom{y_i}{6}} \geq \frac{24}{1555} + o(1).
$$Therefore, as $y_i - y_j \geq 2$ and $x_k = \frac{1}{6} + o(1)$ for all $k \in [6]$, \begin{align*}
N^v + N^{vw} - N^w
	&\leq \frac{N(y_i)}{y_i} + \frac{y_1 \cdots y_6}{y_i} + \frac{y_1 \cdots y_6}{y_iy_j} - \frac{N(y_j)}{y_j} - \frac{y_1 \cdots y_6}{y_j} \\
	&= \frac{y_jN(y_i) - y_iN(y_j)}{y_iy_j} + (y_j-y_i+1)\cdot\frac{y_1 \cdots y_6}{y_iy_j} \\
	&\leq \op{\frac{24}{1555} + o(1)}\op{\frac{1}{y_i}\binom{y_i}{6} - \frac{1}{y_j}\binom{y_j}{6}} + (-2+1)\cdot\frac{y_1 \cdots y_6}{y_iy_j} \\
	&\leq \op{\frac{24}{1555 \cdot 6!} + o(1)}\op{y_i^5 - y_j^5} - \frac{y_1 \cdots y_6}{y_iy_j} \\
	&= \op{\frac{24}{1555 \cdot 6!} + o(1)}(y_i - y_j)(y_i^4+y_i^3y_j+y_i^2y_j^2+y_iy_j^3+y_j^4) - \frac{y_1 \cdots y_6}{y_iy_j} \\
	&\leq \op{\frac{24}{1555 \cdot 6!} + o(1)}\cdot 2\cdot4n^4\op{\frac{1}{6} + o(1)}^4 - n^4\op{\frac{1}{6} + o(1)}^4 \\
	&= n^4\op{\frac{1}{6} + o(1)}^4 \cdot \op{\op{\frac{24\cdot8}{1555 \cdot 6!} + o(1)} - 1} < 0,
\end{align*}a contradiction, proving our claim.
\end{proof}

With this claim, the proof of Theorem \ref{net2} is complete. We now turn our attention back to Theorem \ref{net1}: for $k \geq 1$, the unique graph on $6^k$ vertices which maximizes the number of induced copies of the net graph is a balanced, iterated blow-up of the net. We will note that the following argument is largely recyclable to other graphs, given a proof of corresponding ``Theorem 2." In other words, if one can show that large extremal constructions for a given graph $G$ are (nearly) balanced blow-ups of $G$, then the following argument can be recycled to show that the $|G|^k$-size extremal constructions are unique, and the structure iterates.

\noindent
\begin{proof}
Let $G$ be a minimum counter-example. If $G$ has the outer structure of the net, then the subgraphs inside the blobs are balanced, iterated blow-ups, implying that $G$ would be a balanced, iterated blow-up of the net; as such, $G$ does not have the outer structure of the net. Specifically, $G$ has $6^k > n_0$ vertices, where $n_0$ comes from Theorem \ref{net2}, and $G$ has at least as many nets as the balanced, iterated blow-up $N_{6^k}$.

Take any extremal $H$ on $6^\ell \geq n_0$ vertices, and replace every vertex in $G$ with a copy of $H$ to construct $G_1$, which has at least
$
6^k \cdot N(H) + N(G) \cdot (6^\ell)^6
$ nets. Replace every vertex in $N_{6^k}$ to make $G_2$ and note that it has
$
6^k \cdot N(H) + N(N_{6^k})\cdot(6^\ell)^6,
$ which is extremal by Theorem 2. This implies that $G_1$ is also extremal, so by Theorem 2 we get that $G_1$ is a balanced blow-up of the net with blobs $X_1,\dots,X_6$. Two vertices in $G_1$ are in the same set $X_i$ if and only if their adjacency pattern agrees on more than half of the remaining vertices. But this implies that every copy $H'$ of $H$ in $G_1$ inserted into the blow-up of $G$ has the property that $V(H')$ are coblobular, inducing a net-structure on $G$, a contradiction.
\end{proof}

The only remaining components of the net inducibility problem involve determining all small graphs which maximize the density of induced net graphs. We know that for $n=8$, there is a graph which contains the same number of induced net graphs as a balanced, iterated blow-up of the net graph, indicating that the net graph is not a fractalizer. Therefore, we will not continue to pursue the net inducibility problem beyond these results.

\section{Conclusion}

We note also that Claim \ref{claim5} could be automated with a computer algebra system. 
We have a particular interest in this observation since we suspect that asymmetric graphs may be a fruitful direction of study to find a nontrivial fractalizer, but many simplifying techniques used in the literature are not applicable to them.
A graph is said to be \textit{asymmetric} if its automorphism group is trivial.
This class of graphs is interesting since many counterexamples to being a fractalizer take advantage of graph automorphisms.
For example, it can be shown that any graph with twins is not a fractalizer. The smallest nontrivial assymmetric graphs are on 7 vertices, which is still potentially tractable with flag algebras, and we expect that this will be where the smallest nontrivial fractalizer will be found. Nevertheless, we will not be so bold as to officially conjecture such a thing, and so we instead provide the following weaker conjecture:
\begin{conjecture}\label{ourconj}
There exists an asymmetric fractalizer on at most $9$ vertices. 
\end{conjecture}

We make the above conjecture because (a) we believe it to be true, but also (b) we would certainly like to be proven wrong if possible.

\bibliographystyle{abbrv}
\bibliography{bibliography}

\begin{thebibliography}{10}

\bibitem{sage}
{\em SageMath, the Sage Mathematics Software System (Version 6.3)}, 2020.
\newblock https://www.sagemath.org.

\bibitem{balHuLidPfen}
J.~Balogh, P.~Hu, B.~Lidick\'{y}, and F.~Pfender.
\newblock Maximum density of an induced 5-cycle is achieved by an iterated
  blow-up of a 5-cycle.
\newblock preprint, 2018.

\bibitem{myTournaments}
D.~Burke, B.~Lidick\'{y}, F.~Pfender, and M.~Phillips.
\newblock Inducibility of 4-vertex tournaments.
\newblock Manuscript.

\bibitem{choi}
I.~Choi, B.~Lidick\'{y}, and F.~Pfender.
\newblock Inducibility of directed paths.
\newblock {\em Discrete Mathematics}, 343, 2020.

\bibitem{dary}
E.~Czabarka, A.~A.~V. Dossou-Olory, L.~A. Sz\'{e}kely, and S.~Wagner.
\newblock Inducibility of {$d$}-ary trees.
\newblock {\em Discrete Math.}, 343(2):111671, 15, 2020.

\bibitem{EvenZohar2014}
C.~Even-Zohar and N.~Linial.
\newblock A note on the inducibility of 4-vertex graphs.
\newblock {\em Graphs and Combinatorics}, 31(5):1367--1380, 2014.

\bibitem{MR2988862}
V.~Falgas-Ravry and E.~R. Vaughan.
\newblock Tur\'{a}n {$H$}-densities for 3-graphs.
\newblock {\em Electron. J. Combin.}, 19(3):Paper 40, 26, 2012.

\bibitem{foxhuanglee}
J.~Fox, H.~Huang, and C.~Lee.
\newblock A solution to the inducibility problem for almost all graphs.
\newblock preprint, 2017.

\bibitem{hirst}
J.~Hirst.
\newblock The inducibility of graphs on four vertices.
\newblock {\em Journal of Graph Theory}, 75, 2013.

\bibitem{huVolec}
P.~Hu, B.~Lidick\'{y}, F.~Pfender, and J.~Volec.
\newblock Maximum density of induced oriented $c_4$ is achieved by iterated
  blow-ups.
\newblock Manuscript.

\bibitem{HMNW}
P.~Hu, J.~Ma, S.~Norin, and H.~Wu.
\newblock The inducibility of oriented stars, 2020.

\bibitem{Huang2014}
H.~Huang.
\newblock On the maximum induced density of directed stars and related
  problems.
\newblock {\em SIAM J. Discrete Math.}, 28(1):92--98, 2014.

\bibitem{c5frac}
B.~Lidick\'{y}, C.~Mattes, and F.~Pfender.
\newblock {$C_5$} is almost a fractalizer.
\newblock preprint, 2020.

\bibitem{k2111}
H.~Liu, O.~Pikhurko, M.~Sharifzadeh, and K.~Staden.
\newblock Stability from graph symmetrisation arguments with applications to
  inducibility, 2020.

\bibitem{pipGol}
N.~Pippenger and M.~Golumbic.
\newblock The inducibility of graphs.
\newblock {\em Journal of Combinatorial Theory (B)}, 19:189--203, 1975.

\bibitem{razborov}
A.~Razborov.
\newblock Flag algebras.
\newblock {\em J. Symbolic Logic}, 72:1239--1282, 2007.

\bibitem{bozyk2020inducibility}
Łukasz Bożyk, A.~Grzesik, and B.~Kielak.
\newblock On the inducibility of oriented graphs on four vertices, 2020.

\bibitem{flagmatic}
E.~R. Vaughan.
\newblock Flagmatic 2.0.
\newblock 2012.

\bibitem{yuster}
R.~Yuster.
\newblock On the exact maximum induced density of almost all graphs and their
  inducibility.
\newblock {\em Journal of Combinatorial Theory, Series B}, 136:81--109, 2019.

\end{thebibliography}

\newpage

\section{Appendix A} % Copy from Appendix A

This appendix includes the proof of Claim \ref{claim4}. It was placed here so as to not distract from the more interesting arguments and because the arguments herein are easily reproducible for other graphs.

\noindent
\begin{proof}
To obtain the bounds in Claim \ref{claim4}, we solve four quadratic programs as in the proof of Claim 4 of \cite{balHuLidPfen}. The objectives are to minimize $x_1$, maximize $x_1$, maximize $x_0$, and maximize $f$, respectively. The constraints are $\ds\sum_{i=0}^6 x_i = 1$ and
$$
2 \sum_{1 \leq i < j \leq 6} x_ix_j - 2 \cdot f - a \sum_{i=1}^6 x_i^2 > 0.000149043538
$$in all four cases. By symmetry, bounds for $x_1$ apply also for $x_2,\dots,x_6$. Recall that we chose $a = 4.99$.

Consider the first program described above:
$$
(P) \left\{ \begin{matrix}
	\text{minimize} & & x_1 \\
	\text{subject to} & & \sum_{i=0}^6 x_i = 1 \\
	 & & 2 \sum_{1 \leq i < j \leq 6} x_ix_j - 2 \cdot f - a \sum_{i=1}^6 x_i^2 > 0.000149043538 \\
	 & & x_i \geq 0 \text{ for }i \in \{0,1,\dots,6\}.
\end{matrix}\right.
$$If $(P)$ has a feasible solution $(S)$, then there exists a feasible solution $S'$ of $(P)$ where
\begin{align*}
S'(x_1) = S(x_1), \text{ } S'(f) = 0, \text{ } S'(x_0) = S(x_0), \\
S'(x_2) = \cdots = S'(x_6) = \frac{1}{5}(1 - S(x_1) - S(x_0)).
\end{align*}Since $x_2,\dots,x_6$ appear only in constraints, we only need to check if the second constraint is satisfied. The left hand side of this constraint can be rewritten as
\begin{align*}
2x_1\sum_{2 \leq i < j \leq 6} x_i + 2 \sum_{2 \leq i < j \leq 6} x_ix_j - a \sum_{1 \leq i < j \leq 6} x_i^2 - 2 f \\
	= 2x_1\sum_{2 \leq i < j \leq 6} x_i - \sum_{2 \leq i < j \leq 6} (x_i - x_j)^2 - (a-4) \sum_{1 \leq i < j \leq 6} x_i^2 - a x_1^2 - 2 f.
\end{align*}
Note that the sum of squared differences in the last line is minimized if $x_i = x_j$ for all $i,j \in \{2,\dots,6\}$. The term $x_2^2 + \cdots + x_6^2$, subject to $x_2 + \cdots + x_6$ being a constant, is also minimized if $x_i = x_j$ for all $i,j \in \{2,\dots,6\}$. Since $f \geq 0$, the term $2 f$ is minimized when $f = 0$. Hence, the second constraint in the program is satisfied by $S'$ and we can add the constraints $x_2 = \cdots = x_6$ and $f = 0$ to bound $x_1$. The resulting program $(P')$ is
$$
(P') \left\{ \begin{matrix}
	\text{minimize} & & x_1 \\
	\text{subject to} & & x_0 + x_1 + 5y = 1 \\
	 & & 10x_1y - (4-a) \cdot 5y^2 - a x_1^2 > 0.000149043538 \\
	 & & x_0,x_1,y \geq 0.
\end{matrix}\right.
$$We solve $(P')$ using Lagrange multipliers with the work delegated to Sage \cite{sage}. Finding an upper bound on $x_1$ is done by changing the objective to maximization.

Similarly, we can set $x_1 = \cdots = x_6 = \frac{1}{6}$ to get an upper bound on $f$. We can set $f = 0$ and $x_1 = \cdots = x_6 = (1-x_0)/6$ to get an upper bound on $x_0$. We omit the details as the arguments are similar to above. Sage scripts for solving the all such programs are provided at \url{https://sites.google.com/view/ablumenthal/}.
\end{proof}

\section{Appendix B}

This appendix contains the remaining 5 cases needed to acquire the bounds listed in Table \ref{netBounds} for Claim \ref{claim6}. In each of theses cases, a vertex $x \in X_0$ is moved to a blob and we determine an upper bound on the number of nets containing $x$ and some funky partner $w$. The cases explored involve the locations of both $x$ and $w$, and the upper bound on the number of nets containing this pair vertices is used to provide a lower bound on the normalized funky degree of $x$. As follows are the five remaining cases with some notes about recommended techniques that could be used in future work when studying the inducibility of graphs on $n\geq6$ vertices.

\begin{case} If $x \in X_1$ and $w \in X_4$, then there are at most $$\op{\frac{1}{6}x_0 + \frac{1}{2}d_f(x)\xmax^3 + \frac{1}{2}d_f(x)^2\xmax^2 + \frac{1}{6}d_f(x)^3\xmax}n^4$$ nets containing $x$ and $w$, and so $d_f(x) \geq 0.0610118$.
\end{case}

\noindent
\begin{proof}
As $xw \not\in E(G)$, we know that $w$ is contained in a non-funky $P_4$ in any net containing both $x$ and $w$. We first count nets with a non-funky $P_4$ contained in $X_4$. If $x$ is a triangle vertex, then $w$ is a pendant in the non-funky $P_4$ containing $w$. Further, $p_x \in X_2 \cup X_3 \cup X_4$. This gives at most $n^4 \cdot \op{\frac{d_f(x) \xmax^3}{3!} + \frac{d_f(x) \xmax^2}{2}}$ nets of this type. Otherwise, $x$ is a pendant, so $t_x \in X_4$, giving at most $n^4 \cdot \frac{d_f(x)^3\xmax}{6}$ nets of this type.

Next, we count the number of nets containing $x$ and $w$ with a blob-induced $P_4$. As $w \in X_4$, we know that $w$ is a triangle vertex and hence $x$ is a pendant. Every net containing $x$ and $w$ must then have a blob-induced bull, giving at most $n^4 \cdot \frac{1}{3}d_f(x)\xmax^3$ nets of this type. We thus have a normalized count of $$\op{\frac{1}{6}x_0 + \frac{1}{2}d_f(x)\xmax^3 + \frac{1}{2}d_f(x)^2\xmax^2 + \frac{1}{6}d_f(x)^3\xmax}n^4$$ nets containing $x$ and $w$, and so $d_f(x) \geq 0.0610118$.
\end{proof}

\begin{case} If $x \in X_1$ and $w \in X_5$, then there are at most $$\op{\frac{1}{6}x_0 + \frac{1}{12}\xmax^4 + d_f(x)\xmax^3 + \frac{1}{4}d_f(x)^2\xmax^2} n^4$$ nets containing $x$ and $w$, and so $d_f(x) \geq 0.0433316$.
\end{case}

\noindent
\begin{proof}
We begin by noting that $xw \in E(G)$, and we now outline all of the configurations of nets containing both $x$ and $w$.

Suppose that $x = p_w$. Then any net containing $x$ and $w$ also contains a non-funky $P_4$ which avoids both $x$ and $w$. As such, either this non-funky $P_4$ is contained in $X_5$ or else is blob-induced in $X_1$, $X_4$, $X_6$, and $X_3$. Note that the former case involves no additional funky partners of $x$ and the latter involves precisely one additional funky partner of $x$ located in $X_1$.

Suppose instead that $x = t_w$. Once more, any net containing $x$ and $w$ also contains a non-funky $P_4$ which avoids both of these vertices and must hence be contained entirely in $X_1$ or in $X_3$. Note that the former case involves no additional funky partners of $x$ and the latter involves 2 additional funky partners of $x$ located in either $X_1$ or $X_3$.

Finally, suppose that $x$ and $w$ are both triangle vertices. Note that the third triangle vertex $(t_3)$ must be contained in an inner blob. There is precisely one construction when $t_3 \in X_4$, and it involves precisely one additional funky partner of $x$, its pendant, located in $X_3$. There is precisely one construction when $t_3 \in X_6$, and it involves precisely one funky partner of $x$, in this case $t_3$, located in $X_6$. There are two possible configurations when $t_3 \in X_5$. If $p_x \in X_1$, then we get 1 configuration, and it involves precisely one additional funky partner of $x$, $t_3$, in $X_4$. Otherwise, $p_x \in X_3 \cup X_5$ with 2 potential configuration, both involving an additional funky partner of $x$ in $X_5$ so the non-funky $P_4$ which does not contain $x$ is completely contained in $X_5$. The two configurations differ on whether $p_x$ is in $X_5$ or $X_6$, but this pendant is necessarily an additional funky partner of $x$.

We now bound the number of nets by partitioning the nets over the number of additional funky partners of $x$ and bounding each of these values. There are two configurations where 0 additional funky partners are needed, giving $2 \cdot \binom{\xmax n}{4}$ nets. We now notice that there are precisely four possible configurations involving exactly 1 additional funky partner of $x$, and each of these configurations involve the funky partner in a different blob. This implies that each funky partner of $x$ is in at most $\xmax^3$ nets with $x$ and $w$, giving at most $d_f(x)\xmax^3$ nets containing $x$, $w$ and one more funky partner of $x$.

Finally, we consider how the funky partners of $x$ are distributed among the six blobs more carefully while bounding the number of nets which use two additional funky partners of $x$. One of the configurations involves 2 additional funky partners of $x$ in $X_5$, another involves 2 additional funky partners of $x$ in $X_3$, and the last involves an additional funky partner of $x$ in $X_3$ and another in $X_5$. Therefore, the number of nets containing $x$, $w$, and two more funky partners of $x$ is maximized when all funky partners of $x$ are in $X_3 \cup X_5$, and all pairs of such funky partners can contribute to the count. Therefore, we obtain a reasonable upper-bound of at most $\frac{n^4}{4}d_f(x)^2\xmax^2$.

This gives the claimed enumeration and implies in the usual way that $d_f(x) \geq 0.0433316$. We note that in several cases, we were drastically over-counting by presuming that all choices of 4 vertices in one part would induce the desired structure, namely a $P_4$. This is certainly absurd as the inducibility of the $P_4$ has been shown to be at most 0.204513 in \cite{flagmatic}. Future work on the inducibility of graphs of 6 or more vertices may very well involve this sort of approach, but the improvement of the bounds here does not warrant the added complexity of argument.
\end{proof}

\begin{case} If $x \in X_4$ and $w \in X_1$, then there are at most $$\op{\frac{13}{6} d_f(x) \xmax^3 + \frac{1}{8} d_f(x)\xmax^2 + \frac{1}{6} d_f(x)^3 \xmax}n^4$$ nets containing $x$ and $w$, and so $d_f(x) \geq 0.0322447$.
\end{case}

\noindent
\begin{proof}
In this case, as $xw \not\in E(G)$, we know that $w$ is necessarily in a non-funky $P_4$. We now partition the configurations based on whether or not a non-funky $P_4$ containing $w$ is blob-induced or entirely contained in $X_1$.

There are four cases in which a non-funky $P_4$ containing $w$ is blob-induced: \begin{itemize}
	\item $x$ is a triangle vertex and the $P_4$ that doesn't contain $x$ intersects $X_1$, $X_4$, $X_5$, and $X_2$, with the remaining vertex $(p_x)$ in $X_3$,
	\item $x$ is a triangle vertex and the $P_4$ that doesn't contain $x$ intersects $X_1$, $X_4$, $X_6$, and $X_3$, with the remaining vertex $(p_x)$ in $X_2$,
	\item $x$ is a pendant vertex and there is a blob-induced bull intersecting all blobs but $X_2$, or
	\item $x$ is a pendant vertex and there is a blob-induced bull intersecting all blobs but $X_3$.
\end{itemize}Interestingly enough, each of the four cases involves precisely 1 additional funky partner of $x$, and the inclusion of this funky partner in $X_3$, $X_2$, $X_6$ or $X_5$ (corresponding accordingly to the four cases respectively) completely determines the remaining structure. Again, each funky partner of $x$ contributes at most $\xmax^3$ nets containing $x$ and $w$, so we have at most $n^4d_f(x)\xmax^3$ nets in our count thus far.

Now consider all nets containing $x$ and $w$ wherein the non-funky $P_4$ containing $w$ is entirely within $X_1$. If $x$ is a pendant vertex, then in fact five vertices must be chosen from $X_1$, three of which will be additional funky partners of $x$. This gives at most $\binom{d_f(x)n}{3}\xmax n$ nets. Otherwise, $x$ is a triangle vertex and $w$ is a pendant vertex, so we further partition based on whether $p_x$ is a funky partner of $x$.

If $p_x$ is a funky partner of $x$, then $p_x \in X_2 \cup X_3$. In this case, we need to select a funky partner of $x$ in $X_1$ and another in $X_2 \cup X_3$, as well as two other vertices in $X_1$. We get the largest count in this case when all funky partners of $x$ are in $X_1 \cup X_2 \cup X_3$, and we only induce a net if we choose one funky partner in $X_1$ and one funky partner in $X_2 \cup X_3$. Therefore, there are at most $\frac{1}{2}\binom{d_f(x) n}{2}$ ways to pick two funky partners of $x$ which could possible lead to valid configurations. This gives an upper bound of $\frac{1}{2}\binom{d_f(x) n}{2}\binom{\xmax n}{2}$ nets in this case.

If $p_x$ is not a funky partner of $x$, then $p_x \in X_1 \cup X_5 \cup X_6$. In this case, we get the most nets when all funky neighbors of $x$ are in $X_1$. We necessarily pick 2 vertices in $X_1$, but we pick 3 whenever $p_x$ is in $X_1$. Therefore, we have at most $d_f(x) n \cdot \binom{\xmax n}{3} + 2 \xmax n \cdot d_f(x) n \cdot \binom{\xmax n}{2}$ nets here.

We thus have an upper bound of $$\op{\frac{13}{6} d_f(x) \xmax^3 + \frac{1}{8} d_f(x)\xmax^2 + \frac{1}{6} d_f(x)^3 \xmax}n^4$$ nets containing $x$ and $w$ in this case, implying in the usual way that $d_f(x) \geq 0.0322447,$ as desired.
\end{proof}

\begin{case} If $x \in X_4$ and $w \in X_2$, then there are at most $$\op{\frac{1}{6}x_0 + \frac{1}{2}d_f(x)\xmax^3 + \frac{1}{2}d_f(x)^2\xmax^2 + \frac{1}{6}d_f(x)^3\xmax}n^4$$ nets containing $x$ and $w$, and so $d_f(x) \geq 0.0349529$.
\end{case}

\noindent
\begin{proof}
In this case, we observe that $xw \in E(G)$. We again determine all configurations, then partition the potentially induced nets over the number of additional funky partners of $x$ in the configuration.

First, suppose $xw$ is a triangle-edge. Note that $t_3$ must be placed blob-distance at most 1 from $w$, and that $p_w$ cannot be blob-distance 1 from either $t_3$ or its pendant. The implication of this is that $t_3 \in X_2$, as are $p_w$ and the pendant of $t_3$. We garner 3 configurations when $p_x \in X_1 \cup X_4 \cup X_6$ is not a funky partner of $x$, giving three configurations in which $x$ has precisely 1 funky partner $(t_3)$ in $X_2$. Alternatively, we garner 2 configurations when $p_x \in X_2 \cup X_3$ is a funky partner of $x$, giving a configuration with two funky partners chosen in $X_2$ or a pair of funky partners chosen in $X_2$ and $X_3$.

Next, we consider the brief case where $x = p_w$. In this case, $w$ is a triangle vertex in an outer blob and all such configurations involve a non-funky bull subgraph. As these bull subgraphs are non-funky and contain a triangle vertex in an outer-blob, all five vertices must be coblobular, leading to precisely 1 configuration and no additional funky partners.

Finally, we consider the case where $x = t_w$, and we consider the placement of the non-funky $P_4$ which does not contain $x$. Observe that if the $P_4$ is coblobular in $X_4$, then we have one valid configuration with no additional funky vertices. If the $P_4$ is coblobular in any of $X_1$, $X_2$, $X_3$ or $X_6$, each of these choices corresponds to a valid configuration with 2 additional funky vertices. The final option is that the $P_4$ is blob-induced in $X_1$, $X_4$, $X_6$ and $X_3$, resulting in precisely 1 configuration and 1 additional funky vertex in $X_1$.

We now bound the number of configurations by partitioning over the number of additional funky vertices. Note that there are two configurations with 0 additional funky partner, giving an upper bound of $2 \binom{\xmax n}{4}$. For those configurations containing only 1 additional funky partner, we note that the funky partner must be in $X_1$ or $X_2$. Every funky partner in $X_2$ contributes at most $3\xmax n \cdot \binom{\xmax n}{2}$ nets as $p_x$ could be chosen in 3 different blobs, whereas every funky partner in $X_1$ can contribute at most $(\xmax n)^3$ nets. As such, the worst-case scenario is when all funky neighbors of $x$ are in $X_2$, and each contributes the maximum possible number of nets, giving at most $\frac{3}{2}d_f(x)\xmax^3n^4$ nets containing $x$, $w$, and precisely 1 additional funky partner of $x$.

Lastly, we carefully consider the configurations in which 2 additional funky partners are used. These configurations involve two funky partners chosen in $X_2$ (and two non-funky partners chosen in $X_2$), two funky partners chosen so that one is in $X_2$ and the other is on $X_3$ (and two non-funky partners chosen in $X_2$), or two funky partners chosen in one of $X_1$, $X_2$, $X_3$ or $X_6$ (and two non-funky partners chosen from the same set). The important observation here is gleaned from the parentheticals: every choice of two funky partners of $x$ produces at most $\binom{\xmax n}{2}$ nets, implying that there are at most $\binom{\xmax n}{2}\binom{d_f(x) n}{2}$ nets containing $x$, $w$, and two additional funky neighbors of $x$. This is potentially achievable when all funky partners of $x$ are located in $X_2 \cup X_3$ when we disregard the likelihood of inducing the necessary structures in the blobs.

Our final count gives at most $$\op{\frac{1}{6}x_0 + \frac{1}{2}d_f(x)\xmax^3 + \frac{1}{2}d_f(x)^2\xmax^2 + \frac{1}{6}d_f(x)^3\xmax}n^4$$ nets containing $x$ and $w$, implying that $d_f(x) \geq 0.0349529$ in the usual way.
\end{proof}

\begin{case} If $x \in X_4$ and $w \in X_5$, then there are at most $$\op{\frac{1}{6}x_0 + \frac{1}{2}d_f(x)\xmax^3 + \frac{1}{2}d_f(x)^2\xmax^2 + \frac{1}{6}d_f(x)^3\xmax}n^4$$ nets containing $x$ and $w$, and so $d_f(x) \geq 0.0504913$.
\end{case}

\noindent
\begin{proof}
Note that $xw \not\in E(G)$. First, we consider the configurations which contain $w$ as a triangle vertex. Note that in each of these configurations, $x$ will be a pendant of some other triangle vertex. Further, if a non-funky $P_4$ containing $w$ is blob-induced, then in fact the configuration contains a blob-induced bull. The two possible configurations are when the bull does not intersect $X_3$ and when the bull does not intersection $X_1$. The former case involves 1 additional funky partner of $x$ located in $X_1$, and the latter involves 1 additional funky partner in $X_6$.

Otherwise, the non-funky $P_4$ containing $w$ is entirely within $X_5$. When $x$ is a triangle vertex in such a configuration, there is 1 additional funky partner of $x$ selected in $X_5$, the pendant vertex that is neither $p_x$ nor $w$. We can then pick $p_x$ as a non-funky partner in $X_1 \cup X_5$ or as a funky partner in $X_3$. When $x$ is a pendant vertex, five vertices must be chosen in $X_5$, 3 of which must be additional funky partners of $x$.

Note that every configuration in this case involves at least 1 additional funky partner of $x$. For those configurations involving precisely 1 additional funky partner of $x$, we observe that each choice of funky partner of $x$ can lead to at most 1 valid configuration: \begin{itemize}
	\item if the funky partner is in $X_1$, then the configuration must involve a blob-induced bull which does not intersect $X_3$,
	\item if the funky partner is in $X_6$, then the configuration must involve a blob-induced bull which does not intersect $X_1$, and
	\item if the funky partner is in $X_5$, then $x$ is a triangle vertex, the non-funky $P_4$ is contained entirely in $X_5$, and $p_x \in X_1$.
\end{itemize}Hence, there are at most $d_f(x)\xmax^3n^4$ nets containing $x$, $w$, and exactly one additional funky partner of $x$. There is only one configuration which involves exactly 2 additional funky partners, and the configuration requires a funky partner in each of $X_3$ and $X_5$. For every pair $(w',w'')$ of funky partners of $x$ where $w' \in X_3$ and $w'' \in X_5$ can produce at most $\binom{\xmax^2 n}{2}$ nets, and there are at most $\frac{1}{2}\binom{d_f(x) n}{2}$ such pairs, giving at most $\frac{1}{8}d_f(x)^2\xmax^2n^4$ nets containing $x$, $w$, and two additional funky partners of $x$.

Finally, there is one configuration involving precisely 3 additional funky partners of $x$ and no configurations with 4 additional funky partners, which gives at most $\binom{d_f(x)n}{3}\xmax n$ nets containing $x$, $w$, and at least 3 additional funky partners of $x$ for at most $$\op{\frac{1}{6}x_0 + \frac{1}{2}d_f(x)\xmax^3 + \frac{1}{2}d_f(x)^2\xmax^2 + \frac{1}{6}d_f(x)^3\xmax}n^4$$ nets containing $x$ and $w$. We then have in the usual way that $d_f(x) \geq 0.0504913$ in this case.
\end{proof}

\end{document}